\documentclass[11pt]{article}
\usepackage{latexsym,euscript,amsmath,amssymb,amsbsy,amsfonts,amsthm,amsopn,amstext,amsxtra,amscd}
\usepackage{epsfig}
\usepackage{graphics}
\usepackage{url}

\usepackage[english]{babel}

\theoremstyle{plain}
\newtheorem{theorem}{Theorem}

\newtheorem{proposition}{Proposition}

\theoremstyle{definition}
\newtheorem{example}{Example}

\newcommand{\N}{{\mathbb N}}

\newcommand{\Z}{{\mathbb Z}}

\newcommand{\D}{\text{$\mathcal{D}$}}

\begin{document}

\centerline{\bf On partitions and periodic sequences} \vskip0,5cm
\centerline{Milan Pa\v st\'eka} \vskip1cm

{\bf Abstract.}{\it In the first part we associate a periodic
sequence to a partition and study the connection the distribution of
elements of uniform limit of the sequences. Then some facts of
statistical independence of these limits are proved.} \vskip1cm

\centerline{\bf Introduction} \vskip0,5cm

This paper is inspired by the papers [CIV], [CV], [K] where the
uniformly distributed sequences of partitions of unit interval are
studied. For technical reasons we will study the systems of finite
sequences of unit interval - equivalent of partitions.

Let $V_N=\{v_N(1), ..., v_N(B_N)\}, N=1,2,3, \dots$, $\lim_{N\to
\infty}B_N=\infty$ be a system of finite sequences of elements of
interval $[0,1]$. We say that this system is {\it uniformly
distributed} if and only if
$$
\lim_{N\to \infty}\frac{1}{B_N} \sum_{n=1}^{B_N} f(v_N(n))= \int_0^1
f(x)dx.
$$
By the standard way we can prove
\begin{theorem}
\label{WKPART} The system $V_N, N=1,2,3, \dots$ is uniformly
distributed if and only if
$$
\lim_{N\to \infty}\frac{1}{B_N} |\{n \le B_N; v_N(n) <x\}| = x
$$
for every $x \in [0,1]$.
\end{theorem}
Notion of uniformly distributed sequence introduced firstly Hermann
Weyl in his famous paper [WEY]. Equivalent of this definition is
that a sequence $v$ of elements of interval $[0,1]$ is {\it
uniformly distributed} if and only if the system of finite sequences
$\{v(1), \dots, v(N)\}, N =1,2,3, \dots$ is uniformly distributed.
This means that the sequence $\{v(n)\}$ is uniformly distributed if
and only if the equality
$$
\lim_{N\to \infty}\frac{1}{N} |\{n \le N; v(n) <x\}| = x
$$
holds for arbitrary $x \in [0,1]$. This can be formulated in
equivalent form: Denote by $\D$ the system of sets $A \subset \N$
that the limit $d(A):=\lim_{N \to \infty}\frac{1}{N}|\{n \le N; n
\in A\}|$ exists. In this case we say that the set $A$ has
asymptotic density and the value $d(A)$ we call the {\it asymptotic
density} of $A$. For the details we refer to [PAS].  Clearly the
sequence $\{v(n)\}$ is uniformly distributed if and only if $\{n \in
\N; v(n) <x\} \in \D$ and $d(\{n \in \N; v(n) <x\})=x$ for each $x
\in [0,1]$.

In 1946 constructed R. C. Buck in the paper [BUC] measure density.
He started from the asymptotic density and applied some methods of
measure theory. Denote $r+(m)=\{n \in \N; n \equiv r \pmod{m}\}$ -
the arithmetic progression,  where $m\in \N, r=0,1,2,\dots$,
$0+(m):=(m)$. If $A \subset \N$ then the value
$$
\mu^\ast(A)= \inf \Big\{\sum_{j=1}^k \frac{1}{m_j}; A \subset
\bigcup_{j=1}^k r_j +(m_j)\Big\}
$$
we shall call {\it Buck's measure density} of $A$. A set $A \subset
\N$ is called {\it Buck measurable} if and only if $\mu^\ast(A) +
\mu^\ast(\N \setminus A)=1$. Denote $\D_\mu$ the system of all Buck
measurable sets. Then $\D_\mu$ is an algebra of sets and $\mu =
\mu^\ast|_{\D_\mu}$ is a finitely additive probability measure on
$\D_\mu$. For the details we refer to [PAS]. A sequence of elements
of $[0,1]$ is called {\it Buck uniformly distributed} if and only if
$\{n \in \N; v(n) <x\} \in \D_\mu$ and $\mu(\{n \in \N; v(n) <x\}) =
x$ for all $x \in [0,1]$. (See also [PAS2].)

We say that a sequence $v$ {\it polyadicly continuous on the set} $A
\subset \N$ if and only if for each $\varepsilon>0$ such positive
integer $m$ exists that for all $n_1,n_2 \in A$ we have
$$
n_1 \equiv n_2 \pmod{m} \Longrightarrow |v(n_1)-v(n_2)| <
\varepsilon.
$$

We say that $v$ is polyadicly continuous if and only it polyadicly
continuous on $\N$.

Denote for a sequence $v$ and $N$ - positive integer
$$
E_N(v)=\frac{1}{N}\sum_{n=1}^N v(n).
$$
If there exists the proper limit
$$
\lim_{N \to \infty} E_N(v):= E(v)
$$
the we say that $v$ has  mean value and $E(v)$ is called the {\it
mean value } of $v$. (In literature is this known also as $(C,1)$
summable  sequence and $(C,1)$ limit.)

We say that a sequence of positive integers $\{k_n\}$ is uniformly
distributed in $\Z$ if and only if for every $m \in \N$ and $r=0,1,
, \dots$ the set $\{n \in \N; k_n \in r+(m)\}$ belongs to $\D$ and
its asymptotic density is $\frac{1}{m}$. (Firstly introduced and
studied by I. Niven in the paper [N]).

\begin{theorem}
\label{WKBUD} A sequence $v$ of elements $[0,1]$ is Buck uniformly
distributed if and only if for each sequence $\{k_n\}$ of positive
integers uniformly distributed in $\Z$ the equality
$$
\lim_{N\to \infty} E_N(g(v(k_n))= \int_0^1 g(x)dx
$$
holds for every continuous function $g$ defined on $[0,1]$.
\end{theorem}
For the proof we refer to [PAS] page 122. \vskip1cm \centerline{\bf
Almost polyadic continuity and } \centerline{\bf almost uniform
convergence.} \vskip0,5cm

The following two propositions proved in [PAS] pages 106, 107, 108
will be useful for the next:

\begin{proposition} Let $v$ be a periodic sequence with the period
$B \in \N$. Then $v$ has mean value
$$
E(v)=\sum_{j=1}^B v(j)
$$
and for each sequence $\{k_n\}$ uniformly distributed in $\Z$ we
have
$$
\lim_{N \to \infty} \frac{1}{N}\sum_{n=1}^N v(k_n)=E(v).
$$
\end{proposition}

Every polyadicly continuous sequence can be uniformly approximated
by the periodic functions thus we have

\begin{proposition} Each polyadicly continuous sequence $v$ has mean
value and
$$
\lim_{N \to \infty}\frac{1}{N}\sum_{n=1}^N v(k_n)=E(v)
$$
for each sequence $\{k_n\}$ of positive integers $\{k_n\}$ uniformly
distributed in $\Z$.
\end{proposition}

A sequence $v$ will be called {\it almost polyadicly continuous} if
and only if for every $\delta>0$ there exists a Buck measurable set
$B\subset \N$ that $\mu(B)\le \delta$ and $v$ is polyadicly
continuous on the set $\N \setminus B$.

\begin{theorem}
\label{WKII} If $v$ is a bounded almost polyadicly
continuous sequence then $v$ has mean value and for each sequence of
positive integers $\{k_n\}$ uniformly distributed in $\Z$ we have
\begin{equation}
\label{UDLIM} \lim_{N \to \infty} \frac{1}{N}\sum_{n=1}^N
v(k_n)=E(v).
\end{equation}
\end{theorem}

{\bf Proof.} Let $\varepsilon >0$. Suppose that $v$ is an almost
polydacly continuous sequence and $C>0$ is such contant that
$|v(n)|<C, n \in \N$. Consider $\delta >0$. Then the exists such $m
\in N, r_1,\dots, r_k \in \{0,\dots, m-1\}$ that
$\frac{k}{m}>1-\varepsilon$ and $v$ is polyadicly continuous on the
set $A=\cup_{j=1}^k r_j+(m)$. Define a sequence $\{v_0(n)\}$ by the
following way: $v_0(n)=v(n)$ for $n \in A$ and $v_0(n)=0$ otherwise.
Clearly
$$
|E_N(v)-E_N(v_0)| \le \frac{C}{N}\sum_{n\le N, n \not\in A}1,
N=1,2,3,\dots.
$$
We see immediately that the sequence $\{v_0(n)\}$ is polyadicly
continuous thus it has mean value. Thus from the above inequality we
get that the distance between the upper and lower limit of $E_N(v)$
is smaller than $2C\delta$. For $\delta \to 0^+$ we get that the
$\lim_{N\to \infty} E_N(v)$ exists. The same reasons provide the
equality (\ref{UDLIM}) also. \qed

By combination of Theorem \ref{WKII} and Theorem \ref{WKBUD} we get:

\begin{theorem}
\label{UDBUD} An almost polyadicly continuous sequence of elements
of $[0,1]$ is Buck uniformly distributed if and only if it is
uniformly distributed.
\end{theorem}

We say that a system of sequences $v_N$ converges {\it almost
uniformly } for $N \to \infty$ to sequence $v$ if and only if for
every $\delta>0$ there exists a Buck measurable set $S \subset \N$
that $\mu(S)<0$ and $v_N$ converges uniformly for $N \to \infty$ to
$v$ on the set $\N \setminus S$.

\begin{proposition}
\label{UCONI} If $v_N, N=1,2,3, \dots$ is a system of polyadicly
continuous sequences of elements of $[0,1]$ and it converges almost
uniformly to a sequence $v$ then $v$ is almost polyadicly continuous
and
$$
\lim_{N \to \infty} E(v_N)=E(v).
$$
\end{proposition}

Each finite sequence $V_N, N=1,2,3,\dots$ can be extended to a
periodic sequence $\omega_N$ by the following way
$$
\omega_N(n) = v_N(j) \Longleftrightarrow n \equiv j \pmod{B_N},
j=1,\dots, B_N, n \in \N.
$$

From Theorem \ref{UDBUD} and Proposition \ref{UCONI} we get
immediately

\begin{theorem} Let the system sequences $\omega_N$ converges almost uniformly to a
sequence $\omega$ for $N \to
\infty$,uniformly for $n \in \N$. Then the system of finite
sequences $V_N, N=1,2,3,\dots$ is uniformly distributed if and only
if the sequence $\omega$ is Buck uniformly distributed.
\end{theorem}

The following fact can be proved directly from Cauchy Bolzano
criterion of uniform convergence:

\begin{proposition}
\label{UC} Let $\sum_{N=1}^\infty a_N$ be a convergent series with
positive elements. Suppose that $\{\alpha_N(n)\}, N=1,2,3, \dots$ is
system of sequences that
$$
|\alpha_N(n)-\alpha_{N+1}(n)| \le a_N, n\in \N.
$$
Then $\{\alpha_N(n)\}$ converges uniformly to a suitable sequence
$\{\alpha(n)\}$.
\end{proposition}

\begin{example} Let $\{B_N\}$ be an increasing sequence of positive
integers such that $B_N | B_{N+1}, N=1,2,3,\dots$. Then the series
$\sum_{N=1}^\infty \frac{1}{B_N}$ converges. If we construct a
system of periodic sequences $\{\omega_N(n)\}$, where $\omega_N(n)$
is periodic modulo $B_N, N=1,2,2,\dots$ such that
$$
r \equiv n \pmod{B_N} \Longrightarrow |\omega_N(r)-\omega_{N+1}(n)|
\le \frac{c}{B_N}
$$
for $n,r \in \N$ and $N=1,2,3,\dots$ for some $c>0$ then
$\{\omega_N(n)\}$ converge to a suitable polyadicly  continuous
sequence for $N \to \infty$ uniformly for $n \in \N$.
\end{example}
\vskip1cm \centerline{\bf Statistical independence} \vskip0,5cm
This
notion is introduced in [R]. Let $v_1,\dots,$ $ v_n$ be a sequences
of elements of $[0,1]$. We say that these sequences are {\it
statistically independent} if and only for continuous functions
$g_1,\dots,g_k$ defined on $[0,1]$ and the sequence $u=g(v_1)\dots
g_k(v_k)$ we have
$$
E_N(u)-E_N(g(v_1))\dots E_N(g_k(v_k)) \to 0
$$
for $N \to \infty$.

For sequences almost polyadicly continuous this has the following
more simple form

\begin{proposition} If $v_1,\dots,v_k$  are almost polyadicly continuous  sequences of
elements of
$[0,1]$ and $u$ has the same sense as above then they are statistically independent if
 and only if
$$
E(u)=E(g_1(v_1))\dots E(g_k(v_k))
$$
for every functions $g_1,\dots,g_k$ continuous on $[0,1]$.
\end{proposition}

\begin{proposition}
\label{RELPRIMP} Let $\{v_1(n)\}$ be a periodic sequence with the
period $M_1$ and $\{v_2(n)\}$ be a periodic sequence with the period
$M_2$. If $(M_1,M_2)=1$ then
$$
E(v_1v_2)=E(v_1)E(v_2).
$$
\end{proposition}
{\bf Proof.} From Chinese reminder theorem we get that every $r_1
\in {1,\dots, M_1}$ and $r_2 \in {1,\dots, M_2}$ there exists
uniquely determined $r \in {1,\dots, M_1M_2}$ that $r \equiv
r_1\pmod{M_1}, r \equiv r_2\pmod{M_1}$. Thus
$$
E(v_1v_2)=\frac{1}{M_1M_2}\sum_{r_1, r_2}v_1(r_1)v_2(r_2)=
E(v_1)E(v_2).
$$
\qed

If $v$ is a periodic sequence then for each function $g$ defined on
the set of values of $v$ the sequence $g(v)$ is also periodic. And
so from Proposition \ref{RELPRIMP} we can conclude:

\begin{proposition} If $v_1,\dots, v_n$ are a periodic sequences
with mutually relative prime periods then these sequences are
statistically independent. \end{proposition}

Every continuous function on $[0,1]$ is uniformly continuous. This
implies that if the system of sequences $v_N$, with elements from
$[0,1]$,  converges almost uniformly to a sequence $v$ for $N \to
\infty$ and $g$ is continuous function on $[0,1]$ then $g(v_N)$
converges uniformly to $g(v)$ also. Therefore Proposition
\ref{UCONI} implies:

\begin{proposition}
\label{UCONII} If $v^{(j)}_N, N=1,2,3, \dots$ is such system of
polyadicly continuous sequences, $j=1,\dots, k$ that for each
$N=1,2,3,\dots$ the sequences $v^{(j)}_N, j=1,\dots,k$ are
statistically independent and for each $j=1,\dots,k$ the system of
sequences $v^{(j)}_N, N=1,2,3, \dots$ converges almost uniformly to
a sequence $v^{(j)}$ then the sequences $v^{(j)}, j=1,\dots,k$ are
statistically independent.
\end{proposition}

Suppose that $P$ is some set of primes and $S(P)$ is the semigroup
generated by $P$. We say that a sequence $v$ is {\it polyadicly
continuous with respect to} $P$ if and only if for each
$\varepsilon>0$ there exists $m \in S(P)$ that
$$
n_1\equiv n_2 \pmod{m} \Longrightarrow |v(n_1)-v(n_2)| < \varepsilon
$$
for $n_1,n_2 \in \N$. It is easy to see that in this case the
sequence $v$ is a uniform limit of periodic sequences with periods
belonging to $S(P)$. Thus we obtain

\begin{theorem} Let $P_1,\dots, P_k$ be a disjoint sets of primes.
If $v_1, \dots, v_k$ are such sequences that $v_j$ is polyadicly
continuous with respect to $P_j$ , $j=1,\dots,k$ then the sequences
$v_1, \dots, v_k$ are statistically independent.
\end{theorem}

\begin{example} If $P_1,\dots, P_k$ are a disjoint sets of primes and
$s_j>1$ then  we can for $n \in \N$ define
$$
\alpha_j(n)=\prod_{p|n, p\in P_j}\big(1-\frac{1}{p^{s_j}}\big),
j=1,\dots, k.$$ The condition $s_j>1$ provides that $\alpha_j$ is
polyadicly continuous with respect to $P_j$ for $j=1,\dots, k$. Thus
these sequences are statistically independent.
\end{example}
\begin{example} If $m_1, \dots ,m_k$ are mutually relative primes
positive integers and we define the sequence $v_j(n)$ as van der
Corput sequence with base $m_j$: for $n= a_rm_j^r+\dots + a_1m_j
+a_0, 0\le a_i \le m_j-1$ we put
$$
v_j(n)= \frac{a_r}{m_j^{r+1}}+\dots + \frac{a_1}{m_j^2} +\frac{a_0}{m_j}.
$$
Then the sequences $v_1, \dots, v_k$ are statistically independent.
\end{example}

 \vskip1cm \centerline{\bf
Construction of uniformly distributed partitions} \vskip0,5cm
Suppose now (with loss of generality) that the finite sequences
$V_N, N=1,2,3.\dots$ are increasing. And
\begin{equation}
\label{border} \lim_{N \to \infty}v_N(1)=0.
\end{equation}
 Let $1= j^N_1, ..., j^N_{k_N}$ be an increasing
subsequence of $\{1,\dots, B_N\}$. Then the sequence $V_N$ can be
decomposed into disjoint subsequences
$$
V_N = V(1, N) \cup \dots \cup V(j^N_{k_N}, N)
$$
where $V(k,N)$ is a sequence of consecutive elements of $V_N$ with
smallest element $v_N(j_k)$, $k=1 , \dots ,k_N$. Denote
$$
\ell(k, N) = v_N(j^N_{k+1})-v_N(j^N_k), 1 \le k < k_N, \ell(k_N, N)=
1- v_N(j^N_{k_N}).
$$
Then (\ref{border}) implies
\begin{equation}
\label{interva} \lim_{N \to \infty} \sum_{k=1}^{k_N}\ell(k,N)=1.
\end{equation}
\begin{theorem}
\label{THEOREMI} Assume that
\begin{equation}
\label{UN1} \lim_{N\to \infty} \ell (k, N) =0
\end{equation}
uniformly for $k$ and
\begin{equation}
\label{UN2} \lim_{N\to \infty} \frac{|V(k, N)|}{\ell(k,N)B_N} =1
\end{equation}
uniformly for $k$. Then the system of finite sequences $V_N,
N=1,2,3,\dots$ is uniformly distributed.
\end{theorem}

{\bf Proof.} Let $x \in (0,1)$ be a fixed number. Then the
conditions (\ref{border}) and (\ref{UN1}) provide that there exists
a positive integer $s_N$ that
$$
\sum_{j=1}^{s_N} \ell(j,N) \le x < \sum_{j=1}^{s_N+1} \ell(j,N)
$$
with exception of finite number of $N$.  This yields
\begin{equation}
\label{limx}
 \lim_{N\to \infty} \sum_{j=1}^{s_N} \ell(j,N)=x.
\end{equation}
Clearly
\begin{equation}
\label{UN3} \sum_{j=1}^{s_N} |V(j,N)| \le |\{n \le B_N; v_N(j) \le
x\}| < \sum_{j=1}^{s_N+1} |V(j,N)|.
\end{equation}
From the condition $(\ref{UN2})$ we get that for $\varepsilon > 0$
there exists $N_0$ that for $N\ge N_0$ and $1\le k \le B_N$ we have
$$
(1-\varepsilon)\ell(k, N)B_N \le |V(k,N)| \le (1-\varepsilon)\ell(k,
N)B_N
$$
and by substitution in (\ref{UN3}) we obtain
$$
(1-\varepsilon)\sum_{j=1}^{s_N}\ell(j, N)  \le \frac{|\{n \le B_N;
v_N(j) \le x \}|}{B_N}\le (1+\varepsilon)\sum_{j=1}^{s_N+1} \ell(j,
N).
$$
From this and (\ref{limx}) we can conclude
$$
\lim_{N \to \infty} \frac{|\{n \le B_N; v_N(j) \le x \}|}{B_N} =x.
$$
\qed

Denote
$$
M_N =\max \{|V(k,N)|; k=1,\dots, k_N\}, m_N =\min \{|V(k,N)|;
k=1,\dots, k_N\}$$ and
$$
 L_N =\max \{\ell(k,N); k=1,\dots, k_N\},
\ell_N =\min \{\ell(k,N); k=1,\dots, k_N\}.
$$

\begin{theorem} If
\begin{equation}
\label{LIM0}
 \lim_{N \to \infty} L_N =0
\end{equation}
and
\begin{equation}
\label{LIM1}
 \lim_{N \to \infty} \frac{M_N\ell_N}{m_NL_N} =1
\end{equation}
then the system if finite sequences $V_N, N=1,2,3,\dots$ is
uniformly distributed.
\end{theorem}
{\bf Proof.} We apply Theorem \ref{THEOREMI}. The condition
(\ref{LIM0}) implies that the condition (\ref{UN1}) is fulfilled.

Clearly the inequalities
\begin{equation}
\label{INEQ}
 \frac{m_N}{L_NB_N} \le \frac{|V(k, N)|}{\ell(k,N)B_N}
\le \frac{M_N}{\ell_NB_N}, k=1, \dots ,k_{N_j}
\end{equation}
hold. Thus for proof that the condition (\ref{UN2}) holds it
suffices to prove that $\lim_{N \to \infty}\frac{M_N}{\ell_NB_N} =1$
because the condition (\ref{LIM1}) provides  that the therm on the
left hand side of inequalities (\ref{INEQ}) has the same limit
points as the right hand side therm. Suppose that the right hand
side therm has a limit point smaller than $1$. Then for suitable
$\alpha <1$ and infinite sequence $\{N_j\}$ we have
$$
|V(k, N_j)| \le \alpha\ell(k,N_j)B_{N_j} ,k=1, \dots ,k_{N_j}
$$
and so we get the contradiction
$$
B_{N_j} = \sum_{k=1}^{k_{N_j}} |V(k,N_j)|\le \alpha
\sum_{k=1}^{k_{N_j}}\ell(k,N_j)B_{N_j} \le \alpha B_{N_j}.
$$
If the right hand side therm has a limit point greater than $1$ then
the left hand side therm has the same limit point and so for
suitable $\beta >1$ and an infinite sequence $\{N_j\}$ the
inequalities
$$
|V(k, N_j)| \ge \beta\ell(k,N_j)B_{N_j} ,k=1, \dots ,k_{N_j}
$$
hold. This yields
$$
B_{N_j} = \sum_{k=1}^{k_{N_j}} |V(k,N_j)|\ge \beta
\sum_{k=1}^{k_{N_j}}\ell(k,N_j)B_{N_j}
$$
and so
$$
1 \ge \beta \sum_{k=1}^{k_{N_j}}\ell(k,N_j).
$$
Therefore $\lim_{j \to \infty}\sum_{k=1}^{k_{N_j}}\ell(k,N_j) \le
\frac{1}{\beta} < 1$ - a contradiction with (\ref{interva}). \qed

Faculty of education, University of Trnava, Priemyseln\'a 4, Trnava,
Slovakia.
\end{document}